\documentstyle[12pt]{article}
\newcommand{\supp}{\mathop{\rm supp}\limits}

\newcommand{\Law}{\mathop{\rm Law}\limits}

\begin{document}

 \begin{center}
{\bf Alternativity of Stable Distributions in Linear Topological Spaces,}\\

\vspace{3mm}

{\bf with criterion.}\\

\vspace{5mm}

 $ {\bf E.Ostrovsky^a, \ \ L.Sirota^b } $ \\

\vspace{4mm}

$ ^a $ Corresponding Author. Department of Mathematics and computer science, Bar-Ilan University, 84105, Ramat Gan, Israel.\\

\vspace{4mm}

E - mail:  eugostrovsky@list.ru\\

\vspace{4mm}
3
$ ^b $  Department of Mathematics and computer science. Bar-Ilan University,
84105, Ramat Gan, Israel.\\

E - mail: \ sirota@zahav.net.il\\

 \vspace{4mm}

 Abstract\\

\end{center}

  \vspace{3mm}

   We give in this short report a very simple proof of Zero-One Law for the stable distributions
  in Linear Topological Spaces (LTS).\par

\vspace{4mm}

{\it  Key words and phrases: } Zero-One Law, Linear Topological Spaces, random variables (r.v.), distribution,
stable distribution, weak limit of measures, independence,
concentration functions and inequalities, countable generated  subspaces (CGS), support of measure, compact imbedded
subspace. \par

 \vspace{6mm}

\section{ Introduction. Notations. Definitions. Statement of problem. Short history. }

 \vspace{3mm}

 Let $ L $  with non-trivial topological structure $ \tau $ and measurable Borelian sigma-algebra $  B $
   be Linear Topological Space (LTS) over the ordinary numerical field $  R. $ \par

   {\it  We do not suppose that $ (L, \tau) $ are  agreed in ordinary sense!} In particular,
 the function $ (\lambda, x) \to \lambda \cdot x, \ \lambda \in R, x \in L $   may be
 discontinuous. This is true, for instance, for the famous  Prokhorov-Skorokhod spaces $  D[0,1]. $  \par

 Further, the conjugate space $ L^* $ may be in general case trivial: $ L^* = \{ 0 \};$ but in this case
formulated above result is also trivial. Therefore, we can and will assume
$ \dim L^* \ge 1;$ may be $ \dim L^* = \infty.$ \par

 Let $ \mu  $ be any probabilistic distribution on the Borelian sets $  B; $  for example,
if $ \eta $ be $  L $ valued Borelian measurable random variable  defined on some probability space
with distribution

 $$
 \mu_{\eta} (A) = {\bf P} (\eta \in A),  \ A \in B, \eqno(1.1)
 $$
then $  \mu $ may be equal to $  \mu_{\eta}. $ \par
 Let also $  K = \{  L_0 \}  $ be some {\it class } of subspaces $  L; $ e.g. closed, open, finite or conversely infinite
dimensional, measurable etc. subspaces of the space $  L. $ \par

\vspace{3mm}

{\bf Definition 1.1. }   {\it  We will say that the probabilistic distribution  $ \mu $  is alternative relatively the class $  K = \{  L_0 \} , $
or equally satisfies Zero-One Law relatively this class, if }

$$
\forall L_0 \in K \ \Rightarrow \ \mu(L_0) = 0 \ {\bf or} \ 1. \eqno(1.2)
$$

\vspace{3mm}

{\bf Definition 1.2.}   {\it  The class  $ K  $ of linear subspaces $ L $ is said to be Countable Generated (CG),  write  $ K = CG, $
if for arbitrary $  L_0 \in K = CG $ there exists a countable (or finite) family of continuous linear functionals
$ \{ b_j \}, \ b_j \in L^*  $  for which }

$$
L_0 = \cap_j \ker b_j, \eqno(1.3)
$$
{\it or equally }

$$
x \in L_0 \ \Leftrightarrow \left[ \forall j \ \Rightarrow b_j(x) = 0 \right]. \eqno(1.3)
$$

 \vspace{3mm}

  {\bf Definition 1.3.} {\it  The non-degenerate distribution  of the $ L $  valued r.v. $ S \ \mu = \mu_S  $  is called
  as ordinary stable, if
  there exists a sequence of i.; i.d. r.v. $ \xi_i, \ i = 1,2,\ldots $  with values in this space $  L  $  and a non-random
  positive numerical sequence $  A(n), \ \lim_{n \to \infty} A(n) = \infty $ such that the r.v. $   S  $ is a weak limit (in distribution)
 as $ n \to \infty $ of weighted (normed) sum }

 $$
 S = w - \lim_{n \to \infty} \frac{\xi_1 + \xi_2  + \ldots + \xi_n}{A(n)}. \eqno(1.4)
 $$

  Of course, Gaussian centered distribution is stable. \par

 \vspace{3mm}

{\bf   We will obtain in the second section the Zero-One Law for stable distributions. } \\

 \vspace{3mm}

  Note that this problem has a long history, see, e.g. \cite{Dineen1}-\cite{Zinn1}. Therein was investigated mainly the
 case of Gaussian measure.\par
 As  regards the non-Gaussian stable distribution, see \cite{Dineen1}, \cite{Janssen1}, \cite{Zinn1} and reference therein,
we intend to generalize these results: ease the conditions, reject the assumption of separability of the space $ L, $
find the necessary and sufficient condition for equality $ {\bf P} (S \in L_0) = 1 $ etc.\par

 \vspace{3mm}

\section{Main result: Zero - One Law for the stable distributions in LTS.}

 \vspace{3mm}

{\bf Theorem 2.1.} {\it Let $ (L,\tau) $ be  the Linear Topological Space in our weak sense, $ S $ be the stable distributed r.v.
with values in $ L $ and the  subspace $ L_0 $ of the space $ L $  is Countable Generated. Then}

$$
{\bf P} (S \in L_0) = 0 \ {\bf or } \ {\bf P} (S \in L_0) = 1. \eqno(2.1)
$$
{\bf Proof.} Let $ \{b_j \}  $ be the sequence of continuous linear functionals on the space $  L  $ which generated
the subspace $  L_0. $  Let also $ \{ \xi_i \},  i = 1,2,\ldots $  be i., i.d. r.v.  from the representation  (1.4). \par
 Denote $  \xi = \xi_1. $  Only two cases are possible: \\

 \vspace{3mm}

{\bf A.}

$$
\forall j \Rightarrow {\bf P}(b_j(\xi) = 0) = 1, \eqno(2.2)
$$

{\bf B.}

$$
\exists j_0 \ \Rightarrow {\bf P}(b_{j_0}(\xi) = 0) = q < 1. \eqno(2.3)
$$

 Let us consider first of all the variant {\bf A.} (2.2). Denote

 $$
 S(n) =  \frac{\xi_1 + \xi_2  + \ldots + \xi_n}{A(n)}. \eqno(2.4)
 $$
 As long as the probability   measure $ {\bf P} $ is sigma  additive,

 $$
 {\bf P} \left( \cap_j \cap_i  \{ b_j (\xi_i)  =0 \} \right) = 1. \eqno(2.5)
 $$
  Since the variable $  S(n) $ is linear combination of $ \xi_i, \ i = 1,2,\ldots,n, $

$$
\forall n \ \Rightarrow {\bf P} (S(n) \in L_0) = 1, \  {\bf P} ( \cap_n \{S(n) \in L_0 \} ) = 1,  \eqno(2.6)
$$

$$
{\bf P} (S \in L_0) = 1. \eqno(2.7)
$$

 We consider now the second option {\bf B} (2.3). Let for definiteness

$$
 {\bf P}(b_{1}(\xi) = 0) = q < 1, \eqno(2.8)
$$
then

$$
 {\bf P}(b_{1}(\xi_i) = 0) = q < 1. \eqno(2.9)
$$

 Since the variables $  b_{1}(\xi_i)  $ are independent and identical distributed, we deduce by means of the theory  of
concentration functions, see for example \cite{Petrov1}, chapter 2, p. 45-46:

$$
{\bf P} ( b_1(S(n)) = 0 ) \le \frac{C(\Law(\xi))}{\sqrt{n}} \to 0, \ n \to \infty. \eqno(2.10)
$$

 Therefore $ {\bf P} ( b_1(S) = 0 )   = 0 $ and hence

$$
{\bf P} (S \in L_0) = 0. \eqno(2.11)
$$

\vspace{3mm}

{\bf Remark 2.1. Criterion.} We obtained on the side the {\it criterion } for equality $ {\bf P} (S \in L_0) = 1 $ for any stable
random variable $  S  $ with value in LTS $  L, $ where recall $  L_0 $ is Countable Generated subspace of $ L $ with
generating sequence of linear functionals $   b_j, \ j = 1,,2, \ldots.$  Namely, if

$$
\forall j \Rightarrow {\bf P}(b_j(\xi) = 0) = 1, \eqno(2.12)
$$
then

$$
{\bf P}( S \in L_0 ) = 1 \eqno(2.13)
$$
otherwise

$$
{\bf P}( S \in L_0 ) = 0. \eqno(2.14)
$$

\vspace{3mm}

{\bf Remark 2.2. Another approach.} Suppose the r.v. $ S $ may be represented as a series

$$
S = \sum_{k=1}^{\infty} \zeta_k \ y_k, \eqno(2.15)
$$
where $ \{ \zeta_k \} $  are independent numerical r.v. and $ \{ y_k \} $ are (vectors) from the subspace $  L_0. $
The convergence in (2.15) is understood in the sense of topology $ \tau $ with probability  one.\par
 For instance, the representation (2.15) is true for Gaussian measure in the space of continuous
functions defined on the compact metric space, i.e. Karunen-Loev series.\par
 If (2.15) is true, then evidently

$$
{\bf P} (S \in L_0) = 0 \ {\bf or } \ {\bf P} (S \in L_0) = 1. \eqno(2.16)
$$

\vspace{3mm}

{\bf Remark 2.3.  Non-stable  measures.} The measure $ \nu  $ defined on the Borelian sets $  B  $ of LTS $  L,  $
not necessary to be bounded, is called quasy-stable, if there exists a stable (Probability) measure $  \mu $ for which
the Radon-Nikodym derivative relative $ \nu $ there exists and is positive and finite almost everywhere:

$$
0 < \frac{d \nu}{d \mu} < \infty \ {\bf a.e.} \eqno(2.17)
$$

 Obviously, arbitrary quasy-stable measure is alternative. It is easy to construct the non-alternative
and following non-quasy stable measures; for example, the distribution a random vector $ \vec{\theta} =(\theta_1, \theta_2) $
 on the two-dimensional plane  $  R^2 $ with characteristical function

 $$
 \phi_{\theta_1, \theta_2} (t_1, t_2) \stackrel{def}{=} {\bf E} \exp( i (t_1 \theta_1 + t_2 \theta_2)) =
 0.5 \left[\exp(t_1^2/2) + \exp( t_2^2/2) \right],
 $$
 then

 $$
 {\bf P} (\theta_1 = 0) =  {\bf P} (\theta_2 = 0) = 1/2.
 $$

  \vspace{3mm}

 {\sc  Open question: \ How are all the quasi-stable measures? } The necessary condition is known: they must be
 alternative. This condition is also sufficient in the finite-dimensional case. \par

 \vspace{3mm}

\section{ About supports of Borelian measures in linear topological spaces.}

 \vspace{3mm}

 Let $  (L, ||\cdot|| ) $ be separable Banach space and $ \mu  $ be Probability Borelian measure on $  L. $ It is known,
see   \cite{Ostrovsky2},  \cite{Buldygin1}, \cite{Ostrovsky3} that there exists a compact imbedded into $  L  $
Banach subspace $  L_1 $ such that $  \mu(L_1) = 1. $\par
 This assertion (the existence of compactly embedded subspace with full measure)
 with at the same proof is true when the linear space $  L $ is countable normed.  But the case of
{\it metrizable }   space $  L  $ is to-day not investigated. \par
 Let us consider the following example. Let $  D  $ be the space of  basic functions on the real axis $  R, $
indeed infinite differentiable functions with compact support equipped with usually topology. Let also
$ \psi = \psi(x), \ x \in R $ be some non-zero function from the space $  D  $ with small support:

$$
\supp \psi \subset [-1/4, + 1/4]. \eqno(3.1)
$$

 Introduce the following discrete distributed random variable $  \kappa $ in this space:

$$
{\bf P} \left(\kappa = T_k \psi  \right)  = 2^{-k}, \ k = 1,2, \ldots. \eqno(3.2)
$$
 Here

 $$
 T_k [ \psi](x) \stackrel{def}{=} \psi(x - k);
 $$
evidently, $  T_k [\psi] (\cdot) \in D. $ \par
 Note that the distribution of the r.v. $  \kappa $ has not compact embedded linear subspace of a full measure.
Actually, the arbitrary compact embedded subspace, say $  W  $ must consists on the functions with uniformly bounded
support, say

$$
\exists M = 2,3, \ldots  \forall  g \in W  = W(M)  \ \Rightarrow \supp g \in [0, 2^{M}], \eqno(3.3)
$$
  But for all the values  $  M  $

 $$
 {\bf P} (\kappa \notin W(M) ) \ge \sum_{ k = M+1}^{\infty} 2^{-k} > 0. \eqno(3.4)
 $$
 Thus, $ \supp \mu_{\kappa} $  can not be contained in any compact embedded subspace of full measure.\par

\end{document}